\documentclass[a4paper, aps, amsmath, 10pt, nofootinbib]{revtex4-2}

\usepackage{packages}
\usepackage{commands}

\begin{document}

\title{A Novel Interpolation-Based Method for Solving the One-Dimensional Wave Equation on a Domain with a Moving Boundary}

\author{Michiel Lassuyt}
\email{lassuyt.m@outlook.com}
\thanks{corresponding author}
\affiliation{KU Leuven KULAK, Department of Physics and Astronomy, Etienne Sabbelaan 53, 8500 Kortrijk, Belgium}

\author{Emma Vancayseele}
\email{emma.vancayseele@gmail.com}
\thanks{corresponding author}
\affiliation{KU Leuven KULAK, Department of Physics and Astronomy, Etienne Sabbelaan 53, 8500 Kortrijk, Belgium}

\author{Wouter Deleersnyder}
\email{wouter.deleersnyder@kuleuven.be}
\affiliation{KU Leuven KULAK, Department of Physics and Astronomy, Etienne Sabbelaan 53, 8500 Kortrijk, Belgium}
\author{David Dudal}
\email{david.dudal@kuleuven.be}
\affiliation{KU Leuven KULAK, Department of Physics and Astronomy, Etienne Sabbelaan 53, 8500 Kortrijk, Belgium}
\author{Sebbe Stouten}
\email{sebbe.stouten@kuleuven.be}
\affiliation{KU Leuven KULAK, Department of Physics and Astronomy, Etienne Sabbelaan 53, 8500 Kortrijk, Belgium}
\author{Koen Van Den Abeele}
\email{koen.vandenabeele@kuleuven.be}
\affiliation{KU Leuven KULAK, Department of Physics and Astronomy, Etienne Sabbelaan 53, 8500 Kortrijk, Belgium}

\begin{abstract}
We revisit the problem of solving the one-dimensional wave equation on a domain with moving boundary. In J. Math. Phys. 11, 2679 (1970), Moore introduced an interesting method to do so. As only in rare cases, a closed analytical solution is possible, one must turn to perturbative expansions of Moore's method. We investigate the then made minimal assumption for convergence of the perturbation series, namely that the boundary position should be an analytic function of time. Though, we prove here that the latter requirement is not a sufficient condition for Moore's method to converge. We then introduce a novel numerical approach based on interpolation which also works for fast boundary dynamics. In comparison with other state-of-the-art numerical methods, our method offers greater speed if the wave solution needs to be evaluated at many points in time or space, whilst preserving accuracy. We discuss two variants of our method, either based on a conformal coordinate transformation or on the method of characteristics, together with interpolation.
\end{abstract}

\maketitle
\section*{Keywords}
One-Dimensional Wave equation, Moving boundary problem, Dynamic Domain, Interpolation methods, Numerical methods, Perturbation method of Moore

\section{Introduction}

The wave equation, describing the propagation of waves in various media, is a fundamental model, applied in many domains such as acoustics, electromagnetism, fluid dynamics and quantum mechanics. Its solution under various conditions is crucial for understanding and predicting wave behaviour in numerous applications~\citep{book:malecki_acoustics,book:geradin_mechanical,book:langenberg_ultrasonic,book:fabrizio_electromagnetism,book:jordan_radiation}. Solutions to this equation with time dependent Dirichlet boundary conditions, where the boundaries are fixed in space, are well known and have been extensively studied~\citep{book:strauss_waveequation}.

However, many practical applications involve scenarios where the boundaries of the physical domain are in motion, adding considerable complexity to the problem. One domain where this moving boundary problem arises is in acoustics, describing the propagation and reflection of sound waves in cavities with moving walls and in flowing fluids~\citep{tsuji2013moving,wang2018hybrid,wang2018particle,sun2012immersed}. Pressure waves and stress waves on a moving domain are encountered in biomechanics and material science~\citep{knobloch2015problems}. Moreover, in geophysics, the problem occurs in ocean waves after tectonic activities and wave runup~\citep{gopalakrishnan1983numerical,takase2011}. Another important case is the reflection of electromagnetic waves on moving mirrors~\citep{cooper1980scattering,Moore}, for example in the mechanism for forming wavepackets in lasers~\citep{Petrov_circle,Law_analytical_solution}. Applications of these moving boundaries include the Doppler effect~\citep{christov2017mechanical_doppler} and the dynamical Casimir effect in quantum mechanics~\citep{maclay2004gedanken,Dodonov}. It is shown that motion, due to the effect on the zero-point field energy, may result in the generation of quanta of the electromagnetic field from a vacuum~\citep{Moore, dewitt, fulling_davies}. 

In \citet{Dodonov}, an extensive summary of the state-of-the-art research on the dynamical Casimir effect is provided. More generally, different approaches have been put forward to solve one-dimensional models with moving left or right boundary, being either analytical, quasi-analytical or numerical. Some analytical solutions stem from a conformal coordinate transformation~\citep{Castagnino_transformation}, mapping the original domain onto a fixed domain while preserving the wave equation. By using the transformation, it is possible to describe the Casimir force and the energy density~\citep{cole_schieve_2001,fulling_davies}, although other Hamiltonian methods can be used as described in~\citep{Dodonov}. A Maple package for calculating the energy density is presented in~\citep{alves_software,alves_software_paper}. 
However, most known analytical solutions are limited to some simple boundary dynamics, for example a linear motion in time~\citep{Nicolai_R_linear, Havelock_diffeq, Balazs, gaffour} or certain oscillating movements~\citep{Law_analytical_solution,wegrzyn_2007}. An inverse method, to find a boundary motion that satisfies a predetermined transformation, is presented in~\citep{Castagnino_transformation, Vesnitskii}. Some quasi-analytical perturbation and asymptotic methods to find the conformal transformation were introduced, which are summarised in~\citep{Dodonov}. In particular, the perturbation method introduced by \citet{Moore} can be applied for more general motions, but is inherently limited to slow dynamics. A numeric geometrical method to find the transformation was proposed in~\citep{cole_schieve_1995, ling_and_bo-zang} and adapted in \citep{caro2024} by performing a Richardson extrapolation. Other ways to solve the one-dimensional problem include an expansion over the instantaneous basis~\citep{Dodonov2,schutzhold,Ruser_2005, Ruser_2006}, and the method of characteristics~\citep{Balazs,Petrov_circle,Petrov_torus,dittrich1994instability,cooper1993asymptotic}.

In this paper, we first re-examine Moore's quasi-analytical perturbation method for slow dynamics. \citet{Moore} stated that the convergence properties of the perturbation series are unknown and that a minimal assumption would seem to be that the boundary position is an analytic function of time. We prove that the requirement of being analytic is not a sufficient assumption for Moore's method to converge, and we illustrate this by means of an example. Second, we introduce a novel numerical approach based on interpolation that is also suitable for fast boundary dynamics, which we will refer to as IMR. In comparison with the state-of-the-art numerical method in~\citep{cole_schieve_1995}, the newly introduced method offers greater speed if the wave solution needs to be evaluated at many points in time or space, while having the same accuracy. Third, the proposed interpolation-based method (IMR) can be adapted to work without the conformal coordinate transformation. In that case, the method of characteristics is used together with interpolation. A comparison between the two variants will also be considered.

The paper is organised as follows. In Sec.~\ref{sec:problem_description}, a brief survey of the conformal coordinate transformation technique is provided. Some exact analytical solutions and the solution provided by Moore's perturbation method are summarised in Sec.~\ref{sec:analytical_methods}. We will illustrate the convergence of Moore's approach with two examples, showing the assumption of an analytic boundary position is not sufficient. Sec.~\ref{sec:numerical_methods} contains the details of the state-of-the-art method and the two newly proposed numerical methods: on one hand, interpolation-based methods based on the conformal coordinate transformation (backtracing method and IMR), and on the other hand based on the method of characteristics (IMC), which side-steps the need for a coordinate transformation. In Sec.~\ref{sec:results}, the methods are validated and their performance is compared. Lastly, Sec.~\ref{sec:conclusion} draws some closing remarks.

\section{Problem description}\label{sec:problem_description}
\begin{figure}
\centering
    \includegraphics{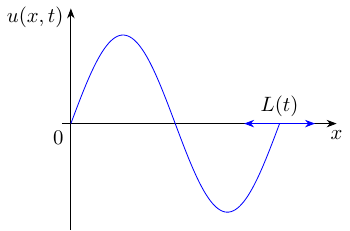}
    \caption{Representation of a wave $u$ on a domain with a variable length.} \label{fig:intro_wave}
\end{figure}

The presented work focuses on the solution of the one-dimensional wave equation on a domain enclosed by a fixed boundary at $x=0$ and a moving boundary at $x=L(t)$, as depicted in Fig.~\ref{fig:intro_wave}. Here, $u(x,t)$ is the displacement at position $x$ and time $t$, and $L(t) \in C^1$ is a function of time $t$ representing the distance between the two boundaries. At the left and right boundaries, we assume homogeneous Dirichlet boundary conditions. We are thus looking for the solution of the system
\begin{align}
    \begin{dcases}
        \pdxx{u}(x,t) = \pdtt{u}(x,t) \qcm & \quad 0 \leq x \leq L(t), \; 0 \leq t \leq t_{\text{max}} \qcm \\
        u(0,t) = 0 \qcm & \quad 0 \leq t \leq t_{\text{max}} \qcm \\[4pt]
        u(L(t),t) = 0 \qcm & \quad 0 \leq t \leq t_{\text{max}} \qcm \\[4pt]
        u(x,0) = f(x)\qcm & \quad 0 \leq x \leq L(0) \qcm \\
        \pdt{u}(x,0) = g(x)\qcm & \quad 0 \leq x \leq L(0) \qcm
    \end{dcases}\label{eq:wave-equation_system}
\end{align}
where the function $f \in C^1$ is the initial distribution of the displacement and $g \in C^0$ is the initial distribution of the velocity. Natural units are being used such that the wave speed equals 1. 

The derivation of the general solution $u(x,t)$ is well-known, and is included here for completeness. For the standard wave equation is given by \[u(x,t) = v(t+x) + w(t-x) \qcm\] where $v$ and $w$ are arbitrary functions. To deal with the variable length boundary condition, \citet{Moore} proposed to map the domain $[0,L(t)]$ onto $[0,1]$ using a conformal coordinate transformation $(x,t) \mapsto (s,\tau)$ of the form
\begin{equation}\begin{cases}
    \tau+s = R(t+x)\qcm\\
    \tau-s = R(t-x) \qpt \label{eq:R(t+-x)}
\end{cases}\end{equation}
The function $R \in C^1$ is chosen to satisfy the boundary condition. This means that the boundary $x=L(t)$ is to be mapped to $s=1$, yielding the following equation that restricts the possible choices for $R$:
\begin{equation}
    R(t+L(t)) - R(t-L(t)) = 2 \qpt \label{eq:R_condition}
\end{equation}
We also demand that $x \in (0,L(t))$ is continuously mapped on $s \in (0,1)$, so $R$ must be strictly increasing. 

Note that the class of transformation~\eqref{eq:R(t+-x)} leaves the wave equation unchanged, since~\citep{Castagnino_transformation} 
\[
\pdxx{} - \pdtt{} =  {R'(t-x)} {R'(t+x)}{\left(\pdss{} -\pdtautau{}\right)} \qpt
\]
As a result, the general solution of the wave equation in $(s,\tau)$-space can be found by separation of variables~\citep{Moore,gaffour}:
\begin{align*}
    u(s,\tau) & = \sum_{\substack{n=-\infty \\ n \neq 0}}^{\infty} C_n \left( e^{-i n\pi (\tau-s)} - e^{-i n\pi (\tau+s)}\right)\qcm
\end{align*}
yielding a solution in $(x,t)$-space given by
\begin{align}
    u(x,t) = \sum_{\substack{n=-\infty \\ n \neq 0}}^{\infty} C_n \left( e^{-i n\pi R(t-x)} - e^{-i n\pi R(t+x)}\right) \qpt \label{eq:u_with_coef_c_n}
\end{align}
The coefficients $C_n$ can be determined given the initial conditions~\citep{Balazs,gaffour}. By inserting these into Eq.~\eqref{eq:u_with_coef_c_n} and extending $f$ and $g$ oddly to $[-L(0),L(0)]$ by setting $f(-x) = -f(x)$ and $g(-x) = -g(x)$, one obtains~\citep{Balazs}
\begin{align}
    C_n = \frac{-i}{4\pi n} \int_{-L(0)}^{L(0)} \left(f'(x)+g(x)\right) e^{i\pi n R(x)} dx \qcm \label{eq:c_n}
\end{align}
such that $C_{-n}$ effectively represents the Fourier series coefficients of the periodic extension of the function 
\begin{align}
\psi(y) = \frac{-i}{\pi n }\frac{f'(R^{-1}(y)) + g(R^{-1}(y))}{2 R'(R^{-1}(y))} \label{eq:fourier_series}
\end{align}
in the interval $y \in [R(-L(0)), 2+R(-L(0))]$. The assumption that the initial conditions are compatible with the boundary conditions means that they must satisfy
\begin{align}
    \begin{cases}
        f(0) = 0 \qcm\\
        g(0) = 0 \qcm\\
        f(L(0)) = 0 \qcm\\
        g(L(0)) = -\dot{L}(0) f'(L(0))\qpt
    \end{cases}\label{eq:initial_cond_pos_velo}
\end{align}
These conditions make sure that the periodic extension of $\psi(y)$, used to calculate the Fourier series, is continuous at the boundary $R(-L(0))$.

Finally, we note that for a simulation from time $t = 0$ to $t = t_{\text{max}}$, it is sufficient to find a transformation $R(\xi)$ on the interval $[-L(0), t_{\text{max}} + L(t_{\text{max}})]$, where we denote $\xi$ as a shorthand for $t \pm x$.

\section{Analytical Methods}\label{sec:analytical_methods}
A (quasi-)analytical solution of system~\eqref{eq:wave-equation_system} can be found for a limited class of boundary dynamics $L(t)$. Straightforward calculations only yield a solution for very specific cases~\citep{Havelock_diffeq, Nicolai_R_linear}, whereas a few more solutions can be found by the inverse method~\citep{Castagnino_transformation, Vesnitskii}. In the case of a slowly moving boundary for which $|\dot{L}(t)| \ll 1$, a perturbation method by \citet{Moore} can be used to approximate the transformation function $R$. 

In this section, we first review the analytical solution for boundary dynamics that can be described by a linear function. In this case it is easy to find the transformation function $R$ explicitly. Next, we shortly recap the inverse method and the perturbation method of Moore which can be applied in more difficult conditions. For the latter method, we illustrate the issue of convergence in two examples: one that converges and one that diverges. In Sec.~\ref{sec:interpol}, we then introduce a novel approach to find the function $R$ using an interpolation-based method.

\subsection{Linear \textit{L}}\label{sec:linear_l}
In the case of a linearly moving boundary $L(t) = L_0 + vt$, it is possible to find an exact solution~\citep{Dodonov, Havelock_diffeq, Nicolai_R_linear}, starting from Eq.~\eqref{eq:R_condition}: $R(t+L_0+vt) - R(t-L_0-vt) = 2$. If $v = 0$, a solution is given by $R(\xi) = \frac{\xi}{L_0}$ resulting in the usual eigenmodes $e^{-i n\pi (t+x)} - e^{-i n\pi (t-x)}$. If $v \neq 0$, we can build the function $R$ by first finding a transformation $\varphi$ such that $t-L_0-vt = \varphi(z)$ and $t+L_0+vt = \varphi(z+1)$. This gives the relation $(1-v)\varphi(z+1) - (1+v)\varphi(z) = 2 L_0$, for which an appropriate choice is given by $\varphi(z) = \frac{-L_0}{v} \pm (\frac{1+v}{1-v})^{z}$. Next, from $(R\circ\varphi)(z+1)-(R\circ\varphi)(z)=2$, it follows that $(R\circ\varphi)(z) = 2z$ up to an arbitrary periodic function of period 1. A solution for $R$ is
\begin{align}
R(\xi) = 2 \frac{\ln\left|1+\frac{v}{L_0}\xi\right|}{\ln\left(\frac{1+v}{1-v}\right)}\qcm \quad \text{if} \; \xi \neq \frac{-L_0}{v} \qpt \label{eq:linear_L}
\end{align}

\subsection{The Inverse Method}\label{sec:inverse_method}
Often it is hard to find a suitable function $R$ corresponding to a given boundary dynamics $L(t)$. Alternatively, one can start from a function $R$ and try to solve Eq.~\eqref{eq:R_condition} inversely in order to extract the underlying function $L$~\citep{Dodonov, Castagnino_transformation, Vesnitskii}. This allows to find simple examples for $L(t)$ which are analytically solvable.

To recover the length $L$, one can in some specific cases algebraically solve Eq.~\eqref{eq:R_condition}. An alternative and more general way is by solving the differential equation 
\begin{equation}
    \dot{L}(t) = \frac{R'[t-L(t)]-R'[t+L(t)]}{R'[t-L(t)]+R'[t+L(t)]} \qcm \label{eq:inverse_method_differential}
\end{equation}
obtained by differentiating Eq.~\eqref{eq:R_condition}, under the constraint $R(L(0)) - R(-L(0)) = 2$. For further reference, the resulting function $L(t)$ for a known transformation $R(\xi)$ is given in 
Tab.~\ref{tab:boundary_conditions}, where the first example is related to the linear boundary dynamics already discussed above, and the second example will be used later for benchmarking. 

\begin{table}
    \caption{Examples of functions $L(t)$ and corresponding $R(\xi)$.}
    \label{tab:boundary_conditions}
    \centering
    \begin{tabular}{ l l l}
    \toprule
    $L(t)$ & \hspace{3mm} & $R(\xi)$ \\
    \cmidrule(lr){1-1}                  
    \cmidrule(lr){3-3}
        $\displaystyle L_0+vt$ && $\displaystyle 2 \frac{\ln\left|1+\frac{v}{L_0}\xi\right|}{\ln\left(\frac{1+v}{1-v}\right)}$ \\[8pt]
        $\displaystyle \frac{1}{k} \arcsinh\left( \frac{1}{A}\sech(k(t-\xi_0)) \right)$ && $\displaystyle A\sinh(k(\xi-\xi_0))$\\[10pt]
       \bottomrule
    \end{tabular}
\end{table}

\subsection{The Perturbation Method of Moore}
\citet{Moore} proposed a general method to find a series equation for $R(\xi)$ in the case of an analytical function $L \in C^{\infty}$ with small boundary velocity $|\dot{L}|$ using perturbation theory. First, the coordinate transformation $s = \frac{x}{L(t)}$ and $\tau = t$ is used to map the time dependent space domain onto a fixed domain. In the case of constant length $L$, it was already reported in Sec.~\ref{sec:linear_l} that $R(t + x)$ is given by $\frac{t + x}{L}$ for each $t\geq0$ and $-L(t) \leq x \leq L(t)$. The time dependent part of this expression is given by $\frac{t}{L}$. If $|\dot{L}(t)| \ll 1$, one still expects that the time dependence of $R$ is mainly characterised by $\int_0^{t} \frac{1}{L(t')}\,dt'$. Let then
\begin{align}
    R(t+x) = R(\tau+sL(\tau)) = g(s,\tau) + \int_0^{\tau} \frac{1}{L(t')}\,dt' \qcm \label{eq:startmoore}
\end{align}
where $g$ is a yet to be specified function of $s$ and $\tau$. In order to find $g$, we remark that the left side of Eq.~\eqref{eq:startmoore} satisfies $\left(\pdt{}-\pdx{}\right)R(t+x) = 0$. Therefore, as a function of the new variables $s$ and $\tau$, this gives a restraint on the function $g$:
\begin{align}
1 + L(\tau) \pdtau{g} = \left( s\dot{L}(\tau) + 1 \right) \pds{g} \qpt\label{eq:g_moore}
\end{align}
Now $g$ can be assumed to be a slowly varying function of time. By defining $\zeta$ as a new time coordinate, only changing significantly in large time intervals: $\zeta = \varepsilon t$, with $\varepsilon \ll 1$ a small constant, and assuming that $g$ can be Taylor expanded in $\varepsilon$ : $g(s,\zeta) = \sum_{\ell=0}^{\infty} g_\ell(s,\zeta)\varepsilon^\ell$, Eq.~\eqref{eq:g_moore} yields
\[
 1 + \varepsilon L(\zeta)\pdzeta{(g_0 + \varepsilon g_1 + \dots)} - \left(s\varepsilon \ddzeta{L(\zeta)}+ 1\right)\pds{\left(g_0 + \varepsilon g_1 + \dots\right)}=0\qpt
\]
By grouping terms with equal powers of $\varepsilon$, a relation between the Taylor series coefficients of $g$, i.e. the functions $g_\ell(s,\zeta)$, is found. Next, each $g_\ell$ is expanded in a power series of $s$ with coefficients $\alpha_{\ell,j}(\zeta)$: $g_\ell = \sum_{j=0}^{\infty} \alpha_{\ell,j} s^j$, which results in a relation between $\alpha_{\ell,j}$ and $\alpha_{\ell,0}$. The boundary condition translates to
\[
g(1,\zeta) - g(-1,\zeta) = 2 \qcm
\]
from which $\alpha_{\ell,0}$ can be found. Finally, by utilizing the Taylor expansion of $\alpha_{\ell,0}$ and rearranging summations, the solution can be expressed in terms of $\gamma_{\ell} = \pdzeta{\alpha_{2\ell-1,0}}$~\citep{Moore} in the original coordinates $\xi=t \pm x$ as follows:
\begin{align}
R(\xi) = \sum^\infty_{\ell=0} \int^\xi_0 \gamma_{\ell}(t') \,dt'\qcm \label{eq:R_moore_sol_final}
\end{align}
where
\begin{align}
    \begin{dcases}
        \gamma_0(t) = \frac{1}{L(t)} \qcm \\
        \gamma_\ell(t) = \sum_{i=1}^\ell \frac{-1}{(2i+1)!} L(t)^{2i} \frac{\mathrm{d}^{2i} \gamma_{\ell-i}}{\mathrm{d}t^{2i}}(t) \quad \text{for } \ell > 0 \qpt
    \end{dcases}
    \label{eq:R_moore_gamma_final}
\end{align}
To calculate $R$ in practice using this method, we need to truncate the series to the first $n$ terms, and we can either calculate the derivatives symbolically or numerically. Computing the derivatives symbolically is expensive and can only be done for simple functions. 
Computing them numerically requires high precision number representation as differentiation is numerically unstable. This means that we can only compute the first few coefficients reliably. As a result, the transformation function $R(\xi)$ we obtain will only be a good fit if $|\dot{L}(t)| \ll 1$, in which case the Taylor series in $\varepsilon$ is expected to converge rapidly.

The true convergence properties of the series are not specifically addressed in previous works~\citep{Moore}. It is simply reasoned that the assumption that $L(t)$ is analytic, would seem a minimal assumption for the convergence, and that it is in practice supposed that the remainder of the truncated series is small~\citep{Moore}. Whereas the remainder is indeed generally small for slow movements, we have found an analytic function $L(t)$ for which the series does not converge, namely $L(t) = e^{-kt}$. Therefore the assumption of \citet{Moore} is not a sufficient condition. Before elaborating on this counter-example, we first illustrate the approximation by Moore for the linear boundary dynamics $L(t)$ in the next subsection.

\subsubsection{Example: linearly varying length}
\begin{figure}
    \begin{subfigure}[t]{.49\linewidth}
        \includegraphics[width = .99\linewidth]{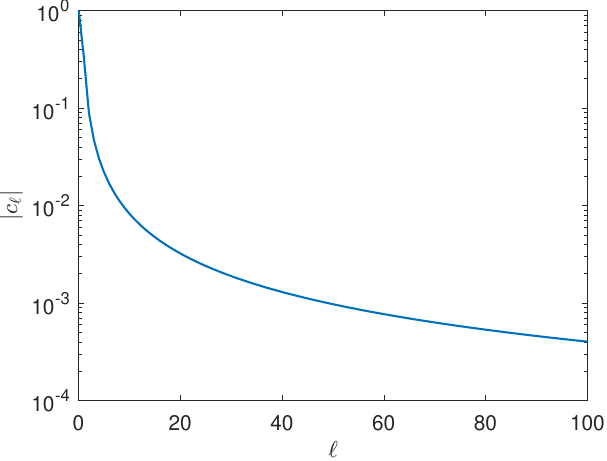}
        \caption{Coefficients $|c_\ell|$ as a function of $\ell$.}
        \label{fig:coefficients_lin}
    \end{subfigure}
    \hfill
    \begin{subfigure}[t]{.49\linewidth}
        \includegraphics[width = .99\linewidth]{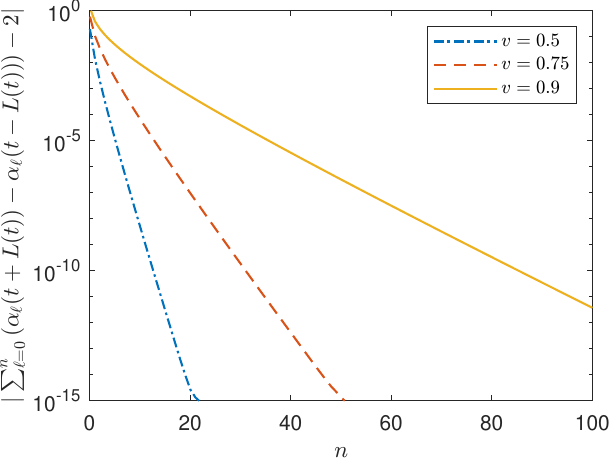}
        \caption{The deviation of the truncated series (Eq.~\eqref{eq:series_coefs}) from the value 2 for different values of the dynamics velocity $v$.}
        \label{fig:convergence_lin}
    \end{subfigure}
    \caption{Convergence of the function $R$ using Moore's method in the case $L(t) = L_0 + vt$.}
\end{figure}

Consider a linearly moving boundary  $L(t) = L_0 + vt$. From the recursion relation \eqref{eq:R_moore_gamma_final}, one can verify by induction that $\gamma_\ell$ has the form $\gamma_\ell = c_\ell v^{2\ell} \frac{1}{L_0 + vt}$. Plugging this into the relation gives
\begin{align*}
\frac{c_\ell v^{2\ell}}{L_0+vt} &= \sum_{i=1}^{\ell} \frac{-1}{(2i+1)!} (L_0+vt)^{2i} \frac{\mathrm{d}^{2i}}{\mathrm{d}t^{2i}} \left(\frac{c_{\ell-i} v^{2\ell-2i}}{L_0+vt} \right) \qcm
\end{align*}
so that the coefficients $c_\ell$ satisfy
\[
c_\ell = -\sum_{i=1}^{\ell} \frac{c_{\ell-i}}{2i+1} \qcm
\]
starting with $c_0 = 1$. Using Eq.~\eqref{eq:R_moore_sol_final}, we can then find $R(\xi) = \sum_{\ell=0}^{\infty} \alpha_\ell(\xi)$ where 
\[
\alpha_\ell(\xi) = \int_0^{\xi}\left( c_\ell v^{2\ell} \frac{1}{L_0+vt}\right)dt = \ln\left(1+\frac{v}{L_0}\xi\right) c_\ell v^{2\ell-1}\qpt
\]
It can indeed be checked that $\sum_{\ell=0}^{\infty} c_\ell x^{2\ell}$ is the unique Taylor series expansion of $\frac{2x}{\ln(\frac{1+x}{1-x})}$ for $|x| < 1$. This means that $R(\xi)$ converges to
\[
R(\xi) = 2\frac{\ln\left(1+ \frac{v}{L_0}\xi\right)}{\ln\left(\frac{1+v}{1-v}\right)}
\]
for $|v| < 1$, which is the same result as in Eq.~\eqref{eq:linear_L}. Figure~\ref{fig:coefficients_lin} shows the decrease of the coefficients $c_\ell$. 
Note that for the boundary condition to be satisfied, we must have $R(t+L_0+vt)-R(t-L_0-vt) = 2$. Denote $R_n$ by the truncated series of $R$, then
\begin{equation}
    R_n(t+L_0+vt)-R_n(t-L_0-vt) = \left[\ln\left(\frac{1+v}{1-v}\right)\right] \sum_{\ell=0}^{n} c_\ell v^{2\ell-1} \label{eq:series_coefs}
\end{equation}
must converge to $2$ for large $n$.
Fig.~\ref{fig:convergence_lin} shows the rate of convergence of this series to the value $2$ as a function of the number of terms $n$ taken into account in Eq.~\eqref{eq:series_coefs}. As expected, the error converges slower for larger speeds $v$, indicating that more coefficients are needed to find an accurate solution. 

\subsubsection{Counterexample: exponentially decreasing length}
\label{sec:moore_coefs_convergence}
\begin{figure}
    \begin{subfigure}[t]{.49\linewidth}
        \includegraphics[width = .99\linewidth]{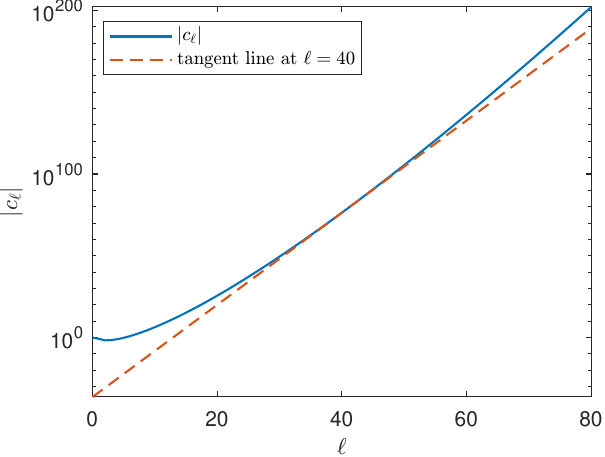}
        \caption{Coefficients $|c_\ell|$ as a function of $\ell$. }
        \label{fig:coefficients_exp}
    \end{subfigure}
    \hfill
    \begin{subfigure}[t]{.49\linewidth}
        \includegraphics[width = .99\linewidth]{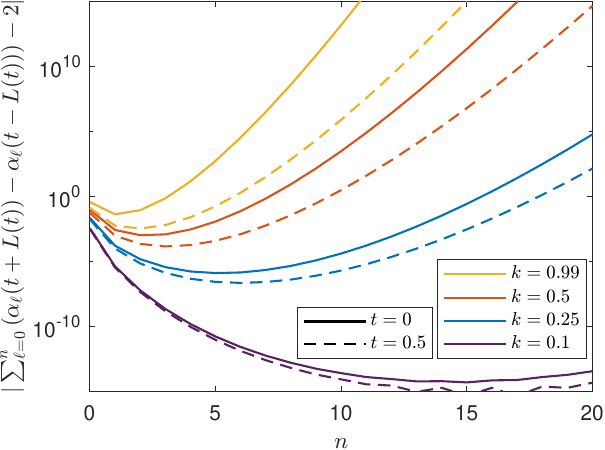}
        \caption{The deviation of the truncated series (Eq.~\eqref{eq:R_alpha_series}) from the value 2 for different values of $k$, evaluated at different time values $t$.}
        \label{fig:convergence_exp}
    \end{subfigure}
    \caption{Convergence of the function $R$ using Moore's method in the case $L(t) = e^{-kt}$.}
\end{figure}

We now analyse Moore's solution for a dynamic boundary condition given by an exponential function $L(t) = e^{-kt}$ where $0<k<1$. This function is analytic and the absolute length change is smaller than 1. In this case, we obtain $\gamma_\ell = c_\ell k^{2\ell} e^{(1-2\ell)kt}$. For the coefficients $c_\ell$, we find
\begin{align*}
c_\ell &= - \sum_{i=1}^{\ell} \frac{(2\ell-2i-1)^{2i}}{(2i+1)!} c_{\ell-i} \qpt
\end{align*}
Fig.~\ref{fig:coefficients_exp} shows that for large $\ell$, the coefficient $c_\ell$ grows in absolute value in the order of ${(a\ell)}^{2\ell}$ with $a > 0$ a constant. In order to quantify this, we first define $c_\ell = (-1)^\ell (2\ell-1)^{2\ell-1} \tilde{c}_\ell$, such that $\tilde{c}_0 = 1$ and for $\ell \geq 1$:
\begin{equation}
\tilde{c}_\ell = \sum_{i=1}^\ell \left(\frac{2\ell-2i-1}{2\ell-1}\right)^{2\ell-1} \frac{(-1)^{i+1}}{(2i+1)!} \tilde{c}_{\ell-i} \qpt \label{eq:c_tilde_expression}
\end{equation}
Next, we will estimate the asymptotic behaviour of $\tilde{c}_\ell$ to determine whether Moore's series converges. For large $\ell$, one easily verifies that only the first few terms in the sum in Eq.~\eqref{eq:c_tilde_expression} are large, in which case $(\frac{2\ell-2i-1}{2\ell-1})^{2\ell-1} \approx e^{-2i}$. In the approximation where we replace this factor in Eq.~\eqref{eq:c_tilde_expression}, it can be shown that the power series  $\sum_{\ell=0}^{\infty} \tilde{c}_\ell (\pi e)^{2\ell} z^{2\ell}$ is equal to $f(z) = \frac{\pi z}{\sin(\pi z)}$ for sufficiently small $z$. Because $f(z)$ is analytic in a disk of radius 1 around 0, and $f(z)$ has a pole at 1, it follows that for large $\ell$, $\tilde{c}_\ell (\pi e)^{2\ell} x^{2\ell}$ converges to zero for $|x|< 1$ but not for $|x| > 1$.
Even without the approximation where the power series equals $f(z)$, we can still obtain roughly the same asymptotic behavior, specifically since $\tilde{c}_\ell (\pi e)^{2\ell} x^{2\ell}$ is not converging to zero for $x > 1$.
Concerning $R$, we find that $R(\xi) = \sum_{\ell=0}^{\infty} \alpha_\ell(\xi)$ where
\begin{align*}
    \alpha_\ell(\xi) = \int_0^{\xi} \left(c_\ell k^{2\ell} e^{(1-2\ell)kt} \right) dt = (-1)^{\ell-1} (2\ell-1)^{2\ell-2} k^{2\ell-1} \tilde{c}_\ell \left( e^{(1-2\ell)k\xi} - 1 \right) \qpt
\end{align*}
Since $(2\ell-1)^{2\ell-2}$ increases faster than $a^{\ell}$ for any $a>0$, the term $\alpha_\ell(\xi)$ does not go to zero, and the infinite series for $R$ diverges for any $0<k<1$ and $\xi \neq 0$. However, as the terms will decrease for small $\ell$ before increasing, only the first few terms of the series can still give a reasonably good approximation. This is similar to the Quantum Electrodynamics perturbation series, which gives excellent comparison with experimental observables using the first few orders, though it is known that the series has convergence radius zero, as first discussed in \cite{dyson1952tj}. For $k$ sufficiently small, which corresponds to a small $|\dot{L}|$, more terms can be added before divergence kicks in, and one can get a reasonably good estimate for $R$. This is seen in Fig.~\ref{fig:convergence_exp}, where the absolute error of $R_n(t+L(t))-R_n(t-L(t))$ with respect to 2 is plotted, as a function of number of terms $n$: 
\begin{equation}
R_n(\xi) = \sum_{\ell=0}^{n} \alpha_\ell(\xi)\qpt \label{eq:R_alpha_series}
\end{equation}

In the next section, we will propose two alternative approaches using numerical methods which do not exhibit this convergence restriction and can effectively handle the case where $|\dot{L}| \mathrel{\not\ll} 1$. A novel interpolation-based method to find the transformation function $R$ and an interpolation-based method of characteristics will be introduced and compared.

\section{Numerical Methods}\label{sec:numerical_methods}
In this section, we aim to reconstruct $R(\xi)$ numerically by exploiting the relation $R(t+L(t)) = 2 + R(t-L(t))$ from Eq.~\eqref{eq:R_condition}. Notice that $t-L(t) \geq -L(0)$ and $t+L(t) \geq L(0)$ for each $t \geq 0$ as $|\dot{L}| < 1$. This means the method will work independent of the values of $R$ on the interval $[-L(0), L(0)]$. The crux of the new method is to iteratively extend $R$ beyond this initial domain. 

For the extension, we can use a polynomial approximation on the initial interval. We demand that $R(-L(0))$ is zero and that $R$, and its derivative, are continuous in $L(0)$. In the most simple case, one could opt for a quadratic polynomial $R(\xi) = \frac{1+\dot{L}(0)}{L(0)}  (\xi + L(0)) -\frac{\dot{L}(0)}{2L(0)^2}  (\xi + L(0))^2$; or, alternatively, a cubic polynomial in order to make the second derivative continuous. More generally, one can also choose a function such that all derivatives in $-L(0)$ and $L(0)$ are zero. This would guarantee a wave solution $u$ that is infinitely differentiable. In this paper, we have performed simulations by employing a cubic polynomial.

\subsection{The Backtracing Method}\label{sec:backtracing}
A first approach, already introduced in literature, e.g. in~\citep{cole_schieve_1995, ling_and_bo-zang}, to iteratively reconstruct $R(\xi)$ outside of $[-L(0), L(0)]$, would be to  state that 
\begin{equation}
    R(\xi) = R(t_1+L(t_1))= R(t_1-L(t_1)) + 2 = R(t_2+L(t_2)) + 2 = R(t_2 - L(t_2)) + 4 \dots \qcm \label{eq:backtrace_procedure}
\end{equation}
using Eq.~\eqref{eq:R_condition}. This process can then be continued until $t_i - L(t_i)$ is part of the initial domain, where the function $R$ can be evaluated. This means that the implicit relation $\xi_j =t_j+L(t_j)$ needs to be solved, which can be challenging in time. 

\subsection{Novel Interpolation-based Methods}
\subsubsection{The Interpolation Method (IMR)}
\label{sec:interpol}

\begin{figure}[b]
\centering
\includegraphics{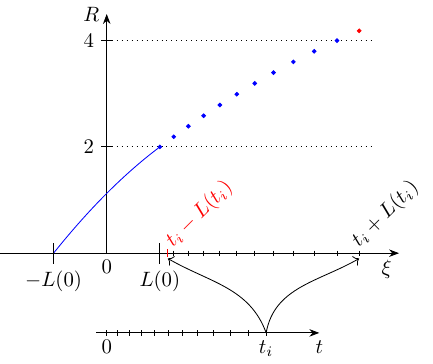}
\caption{Illustration of the interpolation method. Blue indication for $\xi$ values represent instances where $R(\xi)$ is known. New values can be found using the relation \eqref{eq:R_condition}. For example, in order to find the red value at $\xi = t_i + L(t_i)$, it is interpolated at $\xi = t_i - L(t_i)$ based on the blue values, followed by an addition of 2.}
\label{fig:illustration_find_R_interpolation}
\end{figure}

A more elegant approach is by using interpolation to find $R(t_1 - L(t_1))$. We will refer to this method as the Interpolation-based Method to find the transformation function $R$ (IMR). First, a discretisation of the time axis is chosen. For each $t_i$, a corresponding $\xi_i = t_i + L(t_i)$ is defined. By Eq.~\eqref{eq:R_condition}, $R(\xi_i)$ can be calculated from starting from $R(t_i-L(t_i))$, where at first, as illustrated in Fig.~\ref{fig:illustration_find_R_interpolation}, the value $t_i-L(t_i)$ refers to the initial domain for which the function $R$ is known. For later $t_i$, the value may refer to an unknown point outside this domain, and must be estimated. This can be performed by an interpolation on the (discrete) values of $R$ calculated in the previous steps, allowing us to extend the function $R$ to the interval $[-L(0), t_{\text{max}} + L(t_{\text{max}})]$. To avoid extrapolation, we must ensure that $t_i - L(t_i)$ is smaller than the previously calculated $\xi_{i-1} = t_{i-1}+L(t_{i-1})$, in other words $t_i - t_{i-1} < L(t_i) + L(t_{i-1})$. On top of this, the function $R$ generally becomes steeper if $L(t)$ becomes smaller, which implies that we need more points to accurately represent $R$. Thus, instead of employing an equidistant array of $t_i$-values, we choose to work with a time-varying array such that the distance between two consecutive points is proportional to the length of the domain at that time. This way the density of time points is equal to $\rho \frac{1}{L(t)}$, where we will call $\rho$ the resolution. After extending the function $R$, we have acquired an array of function values $R(\xi_i)$ for each discrete value of $\xi_i$. By additional interpolating, we can evaluate the function $R$ in each point $\xi \in [-L(0), t_{\text{max}} + L(t_{\text{max}})]$.

For the interpolation, we use spline interpolation. Notice that the function $R$ consists of several regions that are each a stretched copy of the values on the interval $[-L(0), L(0)]$ shifted vertically by an integer multiple of 2. However, on the boundary between consecutive regions, the function $R$ is not perfectly smooth if we choose to use a polynomial approximation on the initial region: for example, the second derivative is discontinuous in the case of a quadratic polynomial, this will cause a reasonably large interpolation error at these boundaries. 

To resolve this problem, we first iteratively find the time values $\hat{t}$ for which $\hat{t}+L(\hat{t})$ matches the boundary between consecutive regions:
\begin{align*}
\begin{dcases}
    \hat{t}_0 = 0 \\
    \hat{t}_1 \text{ is such that } \hat{t}_1 - L(\hat{t}_1) = \hat{t}_0 + L(\hat{t}_0) = L(0) \\
     \quad\vdots \\
    \hat{t}_{i+1} \text{ is such that } \hat{t}_{i+1} - L(\hat{t}_{i+1}) = \hat{t}_i + L(\hat{t}_i) \qcm
\end{dcases}
\end{align*}
where the corresponding values $\hat{\xi}_i = \hat{t}_i + L(\hat{t}_i)$ lie on the boundaries.
We then add these $\hat{t}$-values to the array of $t$-values assuring that each boundary between regions is represented in the array of $\xi$-values. This allows us to separately use spline interpolation on each region: a spline on the interval $[\hat{\xi}_0, \hat{\xi}_1]$, a spline on the interval $[\hat{\xi}_1, \hat{\xi}_2]$ and so on. This simple adaptation greatly reduces the error of interpolation.

\subsubsection{The Method of Characteristics (IMC)}
\label{sec:method_characteristics}
Instead of determining the transformation function $R$ and calculating the Fourier series coefficients $C_n$ from Eq.~\eqref{eq:c_n} using the initial conditions, we can also directly calculate the wave solution $u(x,t)$, by extending the interpolation-based method introduced in Sec.~\ref{sec:interpol}. We will call this method the Interpolation-based Method of Characteristics (IMC). As a starting point, we write the solution of the wave equation as $u(x,t) = v(t+x) + w(t-x)$. Since the boundary condition $u(0,t) = 0$ implies $w(t) = -v(t)$ for each $t \geq 0$, the wave solution $u(x,t)$ can be rewritten as
\[
u(x,t) = w(t-x) - w(t+x) \qpt
\]
The initial conditions, expressed by $u(x,0) = f(x)$ and $\pdt{u}(x,0) = g(x)$ for each $x \in [0,L(0)]$, imply that $w(-x)-w(x) = f(x)$ and $w'(-x) - w'(x) = g(x)$. If we denote the anti-derivative of $g(x)$ by $G(x)$ with $G(x) = \int_0^x g(x') dx'$, the latter condition requires that $w(-x) + w(x) = -G(x)$ up to a constant. Note that we can simply choose this constant to be zero as it does not affect the solution $u$. Doing so, we obtain the following expressions: $w(x) = \frac{-f(x)-G(x)}{2}$ and $w(-x) = \frac{f(x)-G(x)}{2}$ for each $x \in [0,L(0)]$. Finally, the boundary condition $u(L(t),t) = 0$ results in $w(t+L(t)) = w(t-L(t))$ for each $t\geq 0$. We can now use a similar technique as introduced in Sec.~\ref{sec:interpol} to extend $w$ over a wider domain, thus finding the wave solution $u$ at each time $t \geq 0$.
The crucial difference is that we immediately leverage the initial conditions instead of using the transformation as an intermediate step.

\FloatBarrier
\section{Results and discussion}\label{sec:results}

\begin{figure}
    \centering
    \includegraphics[width=0.44\linewidth]{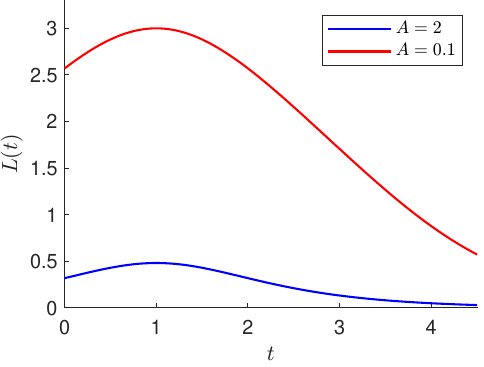}
    \caption{Illustration of the corresponding length evolution $L(t)$ in the case of a hyperbolic sine transformation $R$ as listed in Tab.~\ref{tab:boundary_conditions}, with $k = \xi = 1$. The blue curve is for $A = 2$ and the red curve for $A = 0.1$.} 
    \label{fig:example_sinhR}
\end{figure}

In this section, the objective is to measure the quality of the numerical and analytical methods by means of several error metrics. We will test these methods for two boundary dynamics for which an exact solution of $R$ is available (see Tab.~\ref{tab:boundary_conditions}). The first one corresponds to a linearly moving boundary (discussed in Sec.~\ref{sec:linear_l}) and the second one represents a more complex motion of the boundary involving stretching, squeezing and varying speed. For the second example, the inverse dynamics of a hyperbolic sine function $R$ is chosen, for which $L(t)$ is obtained by means of the inverse method as briefly mentioned in Sec.~\ref{sec:inverse_method}. The length evolution in this case is shown in Fig.~\ref{fig:example_sinhR} for two amplitudes. The exact transformation functions $R$ will be used as a benchmark.

For the initial conditions, we have chosen a gaussian pulse wave and a sine wave, with their respective mathematical expressions summarised in Tab.~\ref{tab:initial_conditions}. The gaussian pulse only interacts with the boundary upon collision, while the standing sine wave is continuously affected by the boundary. The linear initial velocity is chosen as the time derivative of $2\sin(4\pi \frac{x}{L(t)})$ at time zero to satisfy Eq.~\eqref{eq:initial_cond_pos_velo}. Other initial conditions did not significantly impact the simulations (not included for brevity). 

In the next section, we discuss three indicators by which we measure the performance of the different methods. A summary of the compatibility between these indicators and the methods is provided in Tab.~\ref{tab:errors}.

\begin{table}
    \caption{Selected initial conditions used in the comparative study.}
    \label{tab:initial_conditions}
    \centering
    \begin{tabular}{l l l l l}
        \toprule
        Name& \hspace{2mm} &  $f(x)$& \hspace{2mm} & $g(x)$ \\
        \cmidrule(lr){1-1}  \cmidrule(lr){3-3} \cmidrule(lr){5-5}
        Sine wave& &  $ \displaystyle 2\sin\left(4\pi \frac{x}{L(0)}\right)$& & $\displaystyle - \frac{x}{L(0)} \dot{L}(0) f'(x)$\\[12pt]
        Gaussian wave& &$ \displaystyle 2\exp\left(-\tfrac{1}{2} {\left[\left(x-\frac{L(0)}{2}\right) {\mathord{\scalebox{1.5}[1.5]{$/$}}} \frac{L(0)}{16}\right]}^2\right)$& &$\,0$ \\[12pt]
        \bottomrule
    \end{tabular}
\end{table}

\begin{table}
    \caption{Summary of the compatibility of the error indicators and the methods.}
    \label{tab:errors}
    \centering
    \begin{tabular}{l l c c c }
         \toprule
                                        && Moore & IMR & IMC \\
        \cmidrule(lr){3-3} \cmidrule(lr){4-4}  \cmidrule(lr){5-5}
        $\varepsilon_{\text{IC}}$       &Sec.~\ref{sec:results_initial}& $\times$ & $\times$ & \\
        $\varepsilon_{\text{BC},R}$     &Sec.~\ref{sec:results_boundary}& $\times$ & $\times$ & \\
        $\varepsilon_{\text{BC}, w}$    &Sec.~\ref{sec:results_boundary}& & & $\times$ \\
        $\varepsilon_{\text{RMS}}(t)$   &Sec.~\ref{sec:RMSE}& $\times$ & $\times$ & $\times$\\
        \bottomrule
    \end{tabular}
\end{table}

\subsection{Error on the initial conditions}
\label{sec:results_initial}
Given a function $R$ and initial conditions for $t=0$, the wave solution $u$ can be expressed by an infinite linear combination of eigenmodes (see Eq.~\eqref{eq:u_with_coef_c_n}). The coefficients $C_n$ are calculated such that the initial conditions hold. In order to numerically calculate the solution, one must truncate the series and only use the first terms from $n = -n_{\text{max}}$ to $n = n_{\text{max}}$. Thus, even if the right boundary condition is exactly fulfilled at all times, an error can still be present due to the approximation of the initial position and initial speed with a finite linear combination of the eigenmodes. As explained in Appendix~\ref{appendix:error_bound}, a minimal upper bound $\varepsilon_{\text{IC}}$ for this error is given by
\[
\varepsilon_{\text{IC}} = \max_{\xi \in [-L(0),L(0)]} \tilde{w}(\xi) - \min_{\xi \in [-L(0),L(0)]} \tilde{w}(\xi) \qcm
\]
where
\[
\tilde{w}(\xi) = \sum_{\substack{n=-n_{\text{max}} \\ n \neq 0}}^{n_{\text{max}}} C_n e^{-i\pi n R(\xi)} -\frac{1}{2} \left( f(\xi) + \int_0^\xi g(x) dx \right) \qpt
\]

\begin{figure}
    \centering
    \includegraphics[width = .49\linewidth]{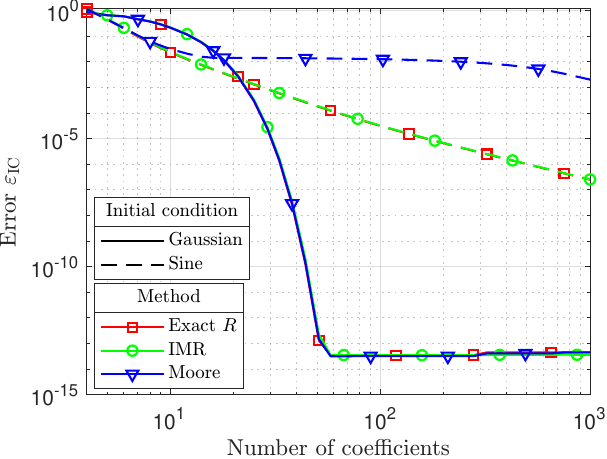}
    \caption{Effect of the number of generalised  Fourier coefficients $C_n$ in Eq.~\eqref{eq:u_with_coef_c_n} on the error bound on amplitude $\varepsilon_{\text{IC}}$, evaluated for three methods: the exact transformation $R$, the newly developed interpolation method (IMR) and Moore's perturbation method; and for two initial conditions: sine and gaussian. The boundary dynamics corresponds to the second example in Tab.~\ref{tab:boundary_conditions} with $A=k=\xi_0=1$, illustrated in Fig.~\ref{fig:example_sinhR} by the blue line.}
    \label{fig:error_bound_on_the_initial_conditions}
\end{figure}

In this first subsection, we compare the effect of the number of terms on the error $\varepsilon_{\text{IC}}$ to determine our best option.
Graphs of the error bounds for different initial conditions can be found in Fig.~\ref{fig:error_bound_on_the_initial_conditions}, calculated for the exact transformation $R$, for the newly introduced interpolation method IMR and for Moore's perturbation method with $n=3$ terms (here we choose the number of terms for which the root mean squared error on the transformation function $R$ is the smallest, see Sec.~\ref{sec:results_boundary}). For the gaussian initial condition, all methods show a substantial decrease in error to the level of the machine precision. This happens for approximations using $10$ to $40$ coefficients, where the eigenmode's wavelength roughly corresponds to the width of the gaussian pulse. For the sine initial condition, the error decreases more slowly and steadily for both the solution based on the exact function $R$ and for the IMR. For Moore's method the error stagnates due to the boundary condition not being adequately satisfied. However, if a smaller boundary speed $|\dot{L}|$ would be considered, the error for Moore's method would also follow the curve for the other two methods. Based on our simulations, we have seen that the slow convergence is caused by the periodic extension of $\psi(y)$ from Eq.~\eqref{eq:fourier_series} not being continuously differentiable at the boundaries (slowing down the convergence), whereas this is the case for the gaussian. 

For other boundary dynamics $L(t)$, we obtain similar error curves for the same initial conditions as the ones shown in Fig.~\ref{fig:error_bound_on_the_initial_conditions} be it shifted to the left or right, depending on the initial length $L(0)$ and the slope of the transformation $R$ in the domain $[-L(0),L(0)]$. From Eqs.~\eqref{eq:c_n} and \eqref{eq:fourier_series}, it follows that more coefficients are required if $R'$ is small, i.e. if the function $R$ is shallow.

\subsection{Error on the boundary condition}
\label{sec:results_boundary}
The value $\varepsilon_{\text{BC}, R} = |R(t+L(t))-R(t-L(t))-2|$ is a measure of how well the boundary condition is fulfilled at a certain time $t$. When the transformation function is known exactly, this value is evidently zero, but an error can occur in case of the IMR and Moore's method. Similarly $\varepsilon_{\text{BC}, w} = |w(t+L(t))-w(t-L(t))|$ is a measure of the error on the boundary conditions for the method of characteristics.

\begin{figure}
    \centering
    \includegraphics[width=.49\linewidth]{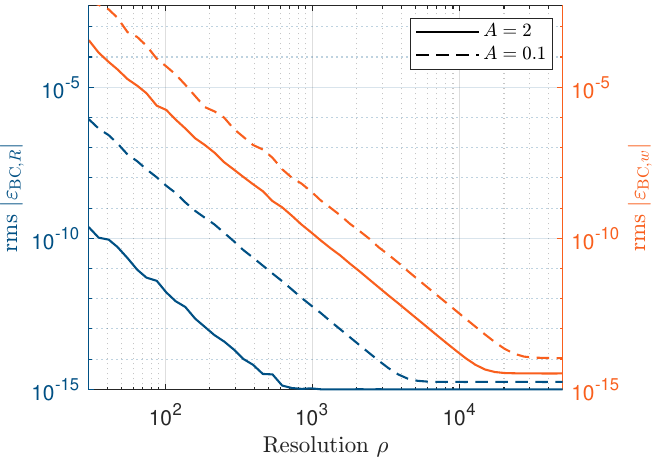}
    \caption{Error on the boundary condition ($\varepsilon_{\text{BC}, R}$ and $\varepsilon_{\text{BC}, w}$) for the IMR and IMC in blue and red respectively, as a function of the resolution $\rho$ (averaged over time using rms). A gaussian was used as initial condition, the boundary dynamics corresponds to the second example in Tab.~\ref{tab:boundary_conditions} with $k=\xi_0=1$, for two values of the parameter $A$. The simulation time $t_{\text{max}}$ is set such that the number of reflections in the simulation are equal in both cases. A larger $A$ value corresponds to a slower moving boundary at lower length. }\label{fig:error_boundary_condition_resolution}
\end{figure}

For the IMR and the IMC, the error $\varepsilon_{\text{BC}}$ is determined by the quality of the interpolation. The influence of the resolution $\rho$ on this error is shown in Fig.~\ref{fig:error_boundary_condition_resolution}. For the IMC, a larger resolution is needed in order to stagnate to machine precision, since the function $R$ is smoother than the function $w$ (determined by the gaussian pulse initial condition). Besides, in order to compare the effect of the parameter $A$, we have set the $t_{\text{max}}$ value such that the number of reflections is equal for both cases, as the error increases with each reflection. We conclude smaller values of the parameter $A$ require a larger resolution due to the larger boundary velocity.

\begin{figure}
    \centering
    \includegraphics[width=.49\linewidth]{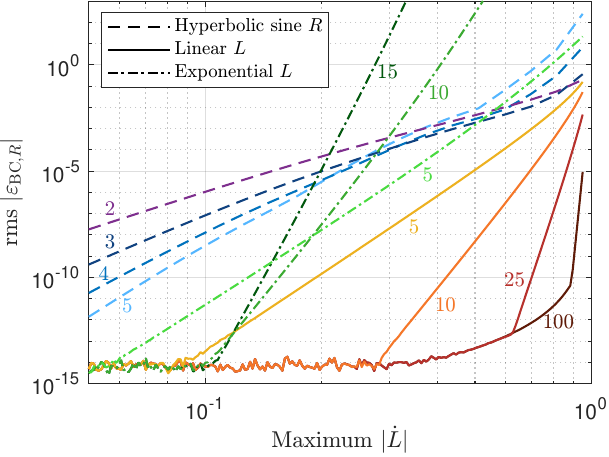}
    \caption{Error on the boundary condition ($\varepsilon_{\text{BC}, R}$)~\eqref{eq:R_condition} for Moore's method as a function of maximal speed of the boundary for different boundary dynamics (averaged over time using rms). The number accompanying the curves specifies the number of terms in the truncation of the series in Eq.~\eqref{eq:R_moore_sol_final}.}
    \label{fig:moore_R_condition}
\end{figure}

For Moore's method, the error is dominated by the convergence of the series in Eq.~\eqref{eq:R_moore_sol_final}, which depends on the velocity of the boundary. Typical error calculations are illustrated in Fig.~\ref{fig:moore_R_condition}. As expected, the error increases for larger $|\dot{L}|$. For a linearly moving boundary, with sufficient number of coefficients, the transformation $R$ satisfies the boundary conditions for both low and high boundary velocities up to machine precision. For other length dynamics, corresponding to a hyperbolic sine $R$ or an exponential $L$, this is not the case. As mentioned in Sec.~\ref{sec:moore_coefs_convergence} and illustrated in Fig.~\ref{fig:moore_R_condition}, there is an optimal number of coefficients, and consequently a minimal achievable error. The latter increases fast for a higher boundary velocity $|\dot{L}|$. Clearly, in these cases other solution methods, such as IMR or IMC, should be considered.

Alternatively, one can validate the exactness of the transformation function $R$ obtained from the boundary dynamics $L(t)$ by employing the inverse method discussed in Sec.~\ref{sec:inverse_method}. Starting from $\hat{R}$, which can be considered as an approximation of the exact transformation function $R$, this results in a corresponding $\hat{L}(t)$, which can be compared to the original $L(t)$. Analysing this `inverse' error on the boundary conditions does not yield new conclusions however.

\subsection{Root Mean Squared Error}\label{sec:RMSE}

For the cases where the exact function $R$ is known (from Tab.~\ref{tab:boundary_conditions}), we can compare the wave solution $u$ of all methods simultaneously to the solution $u_{\text{ref}}$ obtained by using the exact function $R$. Doing so, we can introduce a third error metric which we define as the root mean squared error averaged over the spatial $x$ direction (as done in~\citep{RMSE,Tong_relative_RMSE}):
\begin{align}
\varepsilon_{\text{RMS}}(t) &= \sqrt{\frac{1}{N_x}\sum_{x}\big|u(x,t) - u_{\text{ref}}(x,t) \big|^2} \qcm \label{eq:rmse}
\end{align}
where $N_x$ is the number of $x$ points in $[0, L(t)]$. This indicator measures the deviation of the wave solution $u$ at each time $t$ compared to a reference value. The maximum and mean error were also considered but yield similar results.

Due to the slow convergence of the infinite series of eigenmodes as shown in Fig.~\ref{fig:error_bound_on_the_initial_conditions}, especially in the case of the sine initial condition, we need to take into account both the effect of the error due to the boundary condition (i.e. Eq. \eqref{eq:R_condition}) and the effect due to the approximation of the initial conditions. For the solution based on the exact function $R$ from Tab.~\ref{tab:boundary_conditions}, the boundary condition is exactly satisfied. Furthermore, we impose numerically idealized initial conditions such that there is no error due to their approximations. More specifically, for both the sine and gaussian initial conditions, we calculated the first $300$ coefficients of the eigenmodes. Then, the sum of these $300$ eigenmodes at $t=0$ were used as initial condition instead. Doing so, we have a reliable, exact reference solution.

\begin{figure}
    \begin{subfigure}[t]{0.49\linewidth}
        \includegraphics[width =\linewidth]{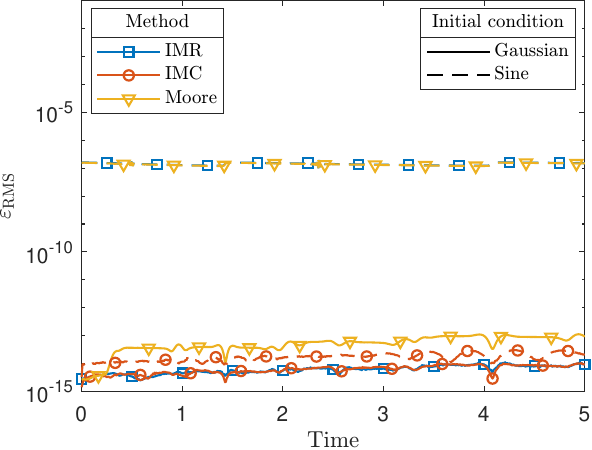}
        \caption{Linearly moving boundary $L(t) = 0.5+0.3t$.}
        \label{fig:RMSE_linear}
    \end{subfigure}
    \hfill
    \begin{subfigure}[t]{0.49\linewidth}
        \includegraphics[width =\linewidth]{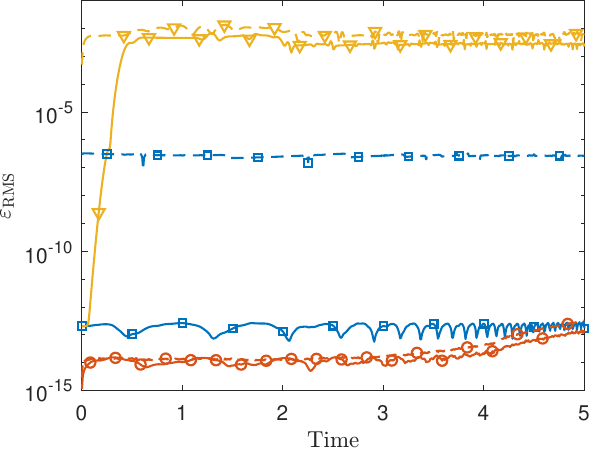}
        \caption{Length corresponding to the hyperbolic sine transformation $R$ with $A = k = \xi_0 = 1$.}
        \label{fig:RMSE_sinh}
    \end{subfigure}
    \caption{Root mean squared error ($\varepsilon_{\text{RMS}}$) compared to the reference solution based on the exact function $R$; for the sine and gaussian initial conditions indicated with full and dotted lines respectively, and for different boundary conditions from Tab~\ref{tab:boundary_conditions}. The methods based on the transformation $R$ were calculated with 150 eigenmodes.}   \label{fig:RMSE}
\end{figure}

Results on the $\varepsilon_{\text{RMS}}(t)$ calculation can be found in Fig.~\ref{fig:RMSE}. The IMC reaches machine precision in all cases, as well as the IMR in the case of the gaussian pulse initial condition. 
For the sine initial condition, the two methods based on the transformation $R$ only reach an error of $10^{-7}$ (in the case of 150 eigenmodes) due to the slow convergence of the coefficients.
Moore's method is comparable with the interpolation method for linearly moving boundary, but performs much worse for the other boundary dynamics: the minimal error only reaches $10^{-2}$ due to a poor compliance of the boundary condition.

The iterative backtracing method developed in~\citep{cole_schieve_1995,ling_and_bo-zang}, as introduced in Sec.~\ref{sec:interpol} and Eq.~\eqref{eq:backtrace_procedure}, is generally quite accurate, as it was found to agree up to machine precision compared to the analytical benchmark. Based on Fig.~\ref{fig:error_boundary_condition_resolution}, we can confirm that the novel IMR is also in agreement up to the machine precision for sufficiently large resolution $\rho$. Hence, our interpolation method is essentially the same as the numerical backtracing method, thus adding the latter to the comparison in this section would not yield new insights.

From this study, we can conclude that the IMC performs the best, followed by the IMR. However, the IMR may be best suited if one is interested in calculating the transformation $R$ (which can be beneficial, for example, for the vacuum energy density calculation which enters the study of the dynamical Casimir effect).

\subsection{Estimate of computation time}\label{sec:interpol_num} 

\begin{figure}
    \includegraphics[width =0.49\linewidth]{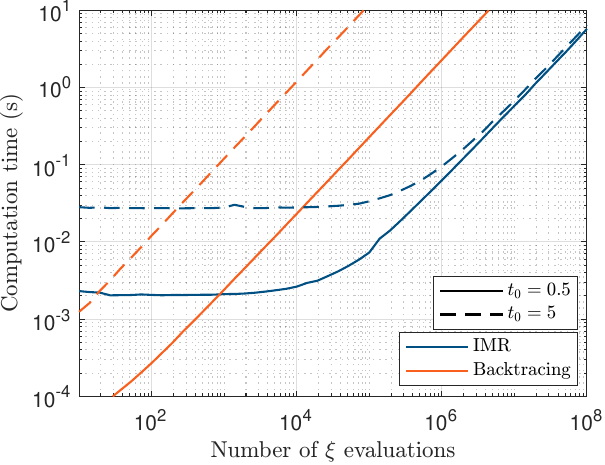}
    \caption{Computation time as a function of the total number of evaluations in space, for two times $t_0$ and for the boundary dynamics corresponding to the second example in Tab.~\ref{tab:boundary_conditions} with $k=\xi_0=1$. An equidistant array in space was used, and for each point $(x_0,t_0)$ the function $R$ was evaluated in $t_0+x_0$ and $t_0-x_0$.}\label{fig:R_comp_time}
\end{figure}

In this final section, we will compare efficiency of the two numerical approaches from Sec.~\ref{sec:interpol}: the iterative backtracing method developed in~\citep{cole_schieve_1995,ling_and_bo-zang} and the newly proposed interpolation method (IMR), both with the coordinate transformation.
We expect the IMR to perform with less computation time than the backtracing method if the transformation function $R$ needs to be calculated at many points. An examination of the computation time as a function of the number of evaluations in space can be found in Fig.~\ref{fig:R_comp_time}, for two times $t_0$. For this graph, an equidistant array in space was used, and for each point $(x_0,t_0)$ the function $R$ was evaluated in $t_0+x_0$ and $t_0-x_0$. The calculations are performed on a personal computer featuring an AMD Ryzen 5 4500U CPU with 2.38 $\si{GHz}$ base clock frequency and 16 $\si{GB}$ of RAM.

Evidently, the computation time increases with the number of evaluations. It is important to note the dominance of the initial cost in computation time for the IMR for a low amount of evaluations, as it needs to find the value of the function $R$ in the interpolation points $\xi_i$ spanning the full range $[-L(0), t_0+L(t_0)]$. Moreover, consideration of a larger $t_0$ value means more preparation time is needed. Eventually, the evaluation phase is linear per evaluation, and not depending on $t_0$. This can be seen at input sizes larger than $10^5$, where the calculation time is solely dominated by the number of evaluations. 

The iterative backtracing method of \citet{cole_schieve_1995, ling_and_bo-zang} does not require a preparation time. The backtracing and evaluation depends on the number of reflections, and increases linearly with input size. For larger $t_0$, more reflections are created, which yields an increase in computation time.

When comparing the two methods, we can conclude that the iterative backtracing method proposed by \citet{cole_schieve_1995, ling_and_bo-zang} is recommended if only a few evaluations are needed. However, many evaluations are often required, for example when simulating a realistic wave where the wave solution $u$ is needed for each $x$ and $t$. In such cases, the IMR is notably faster, with rough estimates showing it to be approximately 25 times faster for $t_0 = 0.5$ and about 5000 times faster for $t = 5$.

\section{Conclusion and outlook}\label{sec:conclusion}
Solving the one-dimensional wave equation under dynamic boundary conditions is a physically most relevant problem. Many analytical methods, as discussed in~\citep{Dodonov}, provide a ground truth. However, their applicability is limited to simple functions. A more general approach is proposed by \citet{Moore} using perturbation theory, though it is only valid for slow dynamics. The present work focused on an interpolation-based method, using a conformal transformation (IMR) or on the method of characteristics (IMC), that can deal with more general scenarios. As a means of validation, we compared the results of the interpolation method to the ones obtained by the method of characteristics for several cases.

For the method proposed by \citet{Moore}, we confirmed that the method performs poorly for rapidly moving boundary dynamics, as it is a perturbation theory with fixed length as unperturbed baseline case, while the other methods successfully deal with these dynamics. Besides, for slower motions, Moore's method performs quite well, but still not as well as the IMC. The here proposed IMR and IMC result definitely in a more accurate solution. Moreover, they also perform well for faster dynamics.   

When comparing the IMR to the iterative method proposed in~\citep{cole_schieve_1995, ling_and_bo-zang}, the trade-off in computation time and accuracy is beneficial for a large amount of evaluations, as the IMR is faster by 1 to 3 orders of magnitude and in agreement with the analytical solution up to machine precision. 

In future work, we believe the IMR and IMC can still be improved. Our choice for the time array 
 suggested in Sec.~\ref{sec:interpol}, namely that the density is proportional to $\frac{1}{L(t)}$, may not be optimal, especially if the length varies near $1$. One can try to find an alternative method to construct the $t$-vector depending on $L$, to minimise the error caused by interpolation. In addition, a different initial transformation function $R$ may be better suited for this purpose.

\section*{CRediT authorship contribution statement}
\textbf{Michiel Lassuyt:} Writing - original draft, Writing - review \& editing, Conceptualisation, Methodology, Software, Validation, Visualization. \textbf{Emma Vancayseele:} Writing - original draft, Writing - review \& editing, Conceptualisation, Methodology, Software, Validation, Visualization. \textbf{Wouter Deleersnyder}: Writing - review \& editing, Supervision. \textbf{David Dudal}: Writing - review \& editing, Supervision. \textbf{Sebbe Stouten}: Writing - review \& editing, Supervision. \textbf{Koen Van Den Abeele}: Writing - review \& editing, Supervision.

\section*{Declaration of competing interest}
The authors declare that they have no known competing financial interests or personal relationships that could have appeared to
influence the work reported in this paper.

\section*{Code availability}
The scripts for the simulations and the implementation of the models are available from \mbox{\url{https://github.com/Bachelor-thesis-KULAK/1D-Moving-Boundary-Problem}}. The programming language is Matlab.

\section*{Acknowledgments}
The work of W.D. was funded by KU Leuven Postdoctoral Mandate PDMt223065. The work of
S.S. was funded by FWO PhD-fellowship fundamental research (file number: 113282). 

\appendix

\section{Error bound on the amplitude}\label{appendix:error_bound}
Suppose the exact solution to the wave equation with the boundary and initial conditions is given by $u_{\text{exact}}$, and the approximate solution using the eigenmodes by $u_{\text{approx}}$. From the wave equation, it follows that $u_{\mathrm{exact}}(x,t) = w_{\mathrm{exact}}(t-x) - w_{\mathrm{exact}}(t+x)$ and $u_{\mathrm{approx}}(x,t) =  w_{\mathrm{approx}}(t-x) - w_{\mathrm{approx}}(t+x)$. We assume now that the only error is caused by the truncated series of the eigenmodes and that we have found an exact $R(\xi)$. Then $\varepsilon =  u_{\text{exact}}-u_{\text{approx}}$ is again a solution of the wave equation with the same boundary conditions, irrespective of the initial conditions, and thus can be written as
\begin{equation}
\varepsilon(x,t) =  w_{\mathrm{exact}}(t-x) - w_{\mathrm{exact}}(t+x) - w_{\mathrm{approx}}(t-x) + w_{\mathrm{approx}}(t+x)\qpt\label{eq:appendix}
\end{equation}
Here $w_{\mathrm{approx}}$ is a finite linear combination of the eigenmodes:
\[
w_{\text{approx}}(\xi) = \sum_{n = -n_{\text{max}}}^{n_{\text{max}}} -C_n e^{-i\pi n R(\xi)} \qpt
\]
After extending $f$ and $g$ to be odd, we obtain from Sec.~\ref{sec:method_characteristics} that 
\[
w_{\text{exact}}(\xi) = \frac{-f(\xi)-G(\xi)}{2} 
\]
for each $\xi \in [-L(0), L(0)]$, where $G(\xi) = \int_0^\xi g(x) dx$. By defining the constant $M$ by
\[
M := \max_{\xi \in [-L(0),L(0)]} \left(w_{\mathrm{exact}}(\xi)- w_{\mathrm{approx}}(\xi)\right) - \min_{\xi \in [-L(0),L(0)]}\left(w_{\mathrm{exact}}(\xi)- w_{\mathrm{approx}}(\xi)\right)\qcm
\]
and using the relation $w(t+L(t)) = w(t-L(t))$, it is clear from Eq.~\eqref{eq:appendix} that the maximal error on the wave $u_{\text{approx}}$ is smaller than $M$.

Next, we define $\xi_1 = \text{argmax}_{\xi \in [-L(0),L(0)]} w(\xi)$ and $\xi_2 = \text{argmin}_{\xi \in [-L(0),L(0)]} w(\xi)$ a,d let $t' = \frac{\xi_1+\xi_2}{2}$ and $x' = \left|\frac{\xi_1-\xi_2}{2}\right|$. Since we assume that $|\dot{L}|<1$, we have for $t \geq 0$ that $t' + x' \leq L(0) \leq t'+L(t')$ and thus $x' \leq L(t')$, which is also true for $t < 0$. This means we can evaluate the wave solution $u$ in $t = t'$ and $x = x'$, and
\[
| u_{\text{exact}}(x',t') - u_{\text{approx}}(x',t') | = M \qcm
\]
which means that the bound $M$ is reached for certain values $x'$ and $t'$, thus it is the smallest bound on the error.

For most functions $R$, this bound will be reached multiple times. Denote $t_0 = 0$ and assume that there exist a $t_1$ such that $t_1-L(t_1) = L(t_0)$. Denote $t_2$ as the time such that $t_2-L(t_2) = L(t_1)$ and so on. Because we assume $L(t)$ to be continuous, the function $w$ passes the same values in the domain $\xi \in [t_i - L(t_i), t_i + L(t_i)]$ as in the domain $\xi \in [-L(0), L(0)]$. So in each interval $[t_i - L(t_i), t_i + L(t_i)]$ there is a $t'$ and an $x'$ such that $| u_{\text{exact}}(x',t') - u_{\text{approx}}(x',t') | = M$.

\bibliographystyle{apsrev4-2}
\bibliography{bibliography}

\end{document}